# Modélisation des données boursières à l'aide de l'analyse de données symboliques

**Philippe Caillou**[*] **— Edwin Diday**[**]

[*] *Ecole Normale Supérieure de Cachan*
61 avenue du Président Wilson, F-94235 Cachan
pcaillou@cybercable.fr

[**] *CEREMADE*
CNRS/Université Paris IX Dauphine
1 place du maréchal de Lattre de Tassigny, F-75016 Paris
diday@ceremade.dauphine.fr

RÉSUMÉ. *Cet article présente une modélisation du domaine boursier à l'aide de l'analyse de données symboliques et une application du logiciel SODAS à ce domaine. Après avoir présenté rapidement le logiciel utilisé, l'analyse est présentée en trois étapes : choix des objets symboliques, définition des objets symboliques et analyse de ces objets à l'aide de SODAS. Chacune de ces étapes est détaillée et leur importance soulignée. Deux exemples de résultats sont ensuite présentés pour montrer l'intérêt et la pertinence de l'analyse. La conclusion décrit les perspectives liées à l'amélioration de SODAS pour l'application au domaine boursier.*

ABSTRACT. *In this paper we present a model of the stock exchange domain using symbolic data analysis and we use the SODAS software to analyze this domain. After a short presentation of the software, we present the analysis in three steps: choice of the symbolic objects, their definition and their analysis with SODAS. We give details for each of these steps and there importance is underlined. Two examples of results are described to show the analysis interest and pertinence. The conclusion describes perspectives after the improvement of SODAS for its application in the stock exchange domain.*

MOTS-CLÉS : *analyse de données symboliques, objet symbolique, bourse, SODAS*

KEY WORDS: *symbolic data analysis, symbolic object, stock exchange, SODAS*

## 1. Introduction

La forte hausse des marchés financiers conjuguée au fort développement d'Internet ont conduit à l'apparition de nombreux sites d'information financières proposant des masses de données considérables, souvent gratuitement. Les données financières et boursières ont de nombreux avantages : elles sont en principe publiques (comptes sociaux des entreprises, cours de bourses, …), elles sont nombreuses et potentiellement à forte valeur ajoutée pour les investisseurs. Elles permettent ainsi aux sites qui les proposent d'attirer un public nombreux, afin éventuellement de proposer des services payants (comme lesechos.fr), ou d'attirer la publicité sur le site (boursorama.com). Ce phénomène a été encore amplifié par le développement des courtiers en lignes, qui proposent des tarifs très faibles pour gérer des comptes et passer des ordres par Internet.

Le particulier ou le professionnel se retrouve ainsi avec une masse de données disponibles sur les différents sites financiers et boursiers telle qu'il ne peut pas les traiter complètement et suffisamment vite. Seul un expert saura où chercher une information pertinente éventuelle car il a une bonne connaissance du domaine et qu'il *s'attend* à trouver l'information là où il va la chercher. Ainsi, un spécialiste du secteur de l'automobile qui sait qu'il y a des mouvements autour du capital de Valeo, société éventuellement opéable, va suivre les volumes de titres échangés, et va effectivement les trouver relativement plus élevés que dans la moyenne du secteur. Une telle information, qui peut être ensuite utilisée par exemple pour anticiper une forte volatilité du titre, est très difficile à rechercher si on ne sait pas à l'avance qu'elle peut se trouver là, car elle est noyée au milieu d'autres données. L'analyse de données, et plus encore l'analyse de données symboliques, apparaissent comme un bon moyen d'extraire des informations pertinentes de ces masses de données brutes.

L'analyse de données peut se définir comme « un ensemble de méthodes dont le but essentiel est la mise en relief des relations existant entre les objets, entre les paramètres qui les caractérisent et entre les objets et les paramètres. » [DID 82]. Il s'agit donc d'un moyen de décrire des objets de façon synthétique à partir et par rapport à des paramètres fournis par l'utilisateur. L'analyse de données symboliques permet d'étendre le champ d'application et les méthodes de l'analyse de données traditionnelle. L'utilisation de données de types quelconques enrichit l'analyse : des données qualitatives, des intervalles, des lois de probabilités ou des variables taxonomiques peuvent ainsi être traitées directement (alors qu'il fallait auparavant les convertir en données numériques standards, perdant ainsi obligatoirement de l'information, ne serait-ce que celle concernant le type de variable). L'analyse de données symboliques met le concept d' « objet symbolique », qui correspond à une description intensionnelle d'une classe d'individus ou d'un individu, au centre de l'analyse. Cela permet d'analyser des classes aussi bien que des individus isolés, mais surtout fournit un cadre formel plus riche que celui d'individu [DID 00][ BOC 82].

SODAS (Symbolic Object Data AnalysiS) est le premier programme d'analyse de données symboliques disponible. Le but de cet article est de montrer comment SODAS et l'analyse de données symboliques peuvent être appliqués aux données boursières et l'intérêt de cette application. On va donc chercher, à l'aide de SODAS, à obtenir une description synthétique et pertinente du marché boursier à partir des données disponibles en grand nombre et gratuitement sur Internet.

## 2. Présentation du logiciel SODAS

SODAS est un logiciel prototype public (accessible à www.cisia.com) apte a analyser des données symboliques. Il est issu du projet de EUROSTAT appelé SODAS comme le logiciel qui en est issu pour fournir un cadre aux différentes avancées récentes et futures du domaine. 11 laboratoires dans 6 pays de la communauté européenne participent à son développement.

Son idée générale est la suivante : à partir d'une base de données, construire un tableau de données symboliques, parfois muni de règles et de taxonomies, dans le but de décrire des concepts résumant un vaste ensemble de données, analyser ensuite ce tableau pour en extraire des connaissances par des méthodes d'analyse de données symboliques. Chaque concept est décrit par des variables dont les valeurs peuvent être des histogrammes, des intervalles, des valeurs uniques (éventuellement munies de règles et de taxonomies) etc., selon le type de variables et le choix de l'utilisateur.

On peut ainsi créer un fichier d'objets symboliques sur lequel une douzaine de méthodes d'analyse de données symboliques peuvent déjà s'appliquer dans le logiciel (histogrammes des variables symboliques, classification automatique, analyse factorielle, analyse discriminante, visualisations graphiques,...).

## 3. Démarche de modélisation

### 3.1. *Choix des objets symboliques*

Cette première étape permet de répondre à la question « Que veut-on décrire ? ». Bien qu'elle paraisse être la plus simple, c'est la plus importante car elle détermine les deux étapes suivantes. L'importance de l'étape dans la démarche globale doit bien être comprise avant d'effectuer le choix des objets qui seront étudiés.

De même qu'en gestion de portefeuille, le choix du marché occupe 10% du temps mais détermine 90% des gains, le choix des objets à décrire en analyse de données peut être rapidement effectué mais aura une grande importance sur la qualité des informations obtenues à la fin. En effet, SODAS ne fait qu'analyser des objets symboliques fournit par l'utilisateur, et fournit une description de cet ensemble d'objets sans chercher à juger de leur qualité.

Que peut-on vouloir décrire en bourse ? L'utilisation des données boursière a généralement pour but le choix de titres pour investir. Dans ce cadre, on peut s'inspirer des méthodes de gestion de portefeuille pour établir une méthodologie d'analyse. La méthode la plus commune en gestion de portefeuille se divise en trois étapes :

- Répartition des fonds entre les marchés (actions françaises, étrangères, obligations, marché monétaire, …). C'est le choix le plus important parce que l'évolution de ces titres financiers est très différente, aussi bien au niveau du rendement que de la volatilité.

- Répartition sur chaque marché entre les secteurs (actions des secteurs informatiques, télécommunications, automobiles, obligations à plus ou moins long terme, …)

- Choix pour chaque secteur de chaque marché des titres à acheter et à vendre. C'est le rôle du trader, la tâche qui prend le plus de temps mais aussi celle qui a l'impact le plus faible, les titres d'un même secteur ayant un comportement voisin.

Pour un débutant comme pour un habitué de la bourse, 4 niveaux d'analyse peuvent donc être intéressants : niveau global, marché, secteur et action. On peut ajouter un cinquième niveau pour quelqu'un qui cherche à avoir des informations sur un ensemble de titres qu'il a déjà ou qu'il a sélectionné : celui du portefeuille. Pour chaque niveau, il faut choisir une granularité qui correspondra à l'objet symbolique. Ainsi, on peut chercher à décrire un marché en terme de temps, de secteur ou d'actions, les objets symboliques considérés seront différents.

Pour chaque niveau d'analyse, il faut donc effectuer un choix de granularité dans les niveaux inférieurs ou par rapport au temps. Chaque combinaison est intéressante et répond à une question différente. Plus la différence de généralité entre objet décrit et granularité est grande, plus le nombre d'objets symboliques sera important. Cela permet de faire ressortir des cas particuliers (tels que des actions au comportement particulier sur la globalité du marché ou un pic soudain sur une longue période) mais rend plus difficile l'interprétation des grands ensembles regroupant un très grand nombre de valeurs. On peut résumer les combinaisons qui semblent les plus intéressantes à l'aide d'un tableau (pour 250 valeurs, ces valeurs sont cotées sur quatre marchés : Règlement mensuel, Règlement mensuel étranger, Second marché et Nouveau marché. Les secteurs repris ici correspondent à la classification officielle en 3 niveaux : 3 secteurs au premier niveau de division, 12 au deuxième et 22 au troisième) :

| Granularité<br>Objet décrit | Marché | Secteur (niveau 3) | action | Temps (semaine) |
|---|---|---|---|---|
| Marché global | comparaison des marchés entre eux<br>*OS : marché (4)* | comparaison inter-sectorielle<br>*OS : secteur (22)* | Description globale par rapport aux titres<br>*OS : action (250)* | Description de l'évolution globale du marché<br>*OS : période (20)* |
| Marché | | Description d'un marché par rapport aux secteurs<br>*OS : secteur (<22)* | Description de la composition d'un marché<br>*OS : action (~100)* | Description de l'évolution d'un marché<br>*OS : période (20)* |
| Portefeuille | description du portefeuille en terme de marché<br>*OS : marché (4)* | Description d'un portefeuille par rapport aux secteurs<br>*OS : secteur (~6)* | Description de la composition du portefeuille<br>*OS : action (~15)* | Description de l'évolution du portefeuille<br>*OS : période (20)* |
| Secteur | | | Description de la composition d'un secteur<br>*OS : action (~15)* | Description de l'évolution d'un secteur<br>*OS : période (20)* |
| Action | | | | Description de l'évolution d'une action<br>*OS : période (20)* |

**Tableau 1.** Présentation des choix d'objets symboliques possibles

**3.2. *Définition des objets symboliques***

Une fois les objets symboliques choisis (action, secteur, …), il faut les définir le plus précisément possible afin que SODAS ait un maximum de données à analyser. Que ce soit pour des objets symboliques marché, secteur ou action, la démarche est ici la même. Un secteur ou un marché reprend les données des actions qui le composent sous forme d'intervalles ou d'attributs modaux multivalués. Aucune donnée spécifique aux secteurs ou aux marchés n'a été considérée.

L'intégration des données dans SODAS s'est faite en quatre étapes : choix de données pertinentes pour décrire une action, recherche de ces données, réunion de celles-ci sous une forme utilisable à l'aide d'Excel et d'Access, et enfin importation de ces données au format SODAS grâce à DB2SO.

3.2.1.*Choix des données*

Les données choisies pour décrire les objets symboliques doivent être à la fois précises et nombreuses pour que l'analyse soit efficace, et pertinentes pour que les informations extraites le soient également. Le programme ne peut en effet pas savoir quelles variables sont pertinentes ou non, c'est à l'utilisateur de lui fournir des données qui le soient. C'est d'autant plus important que la réponse proposée par le programme utilisera peut-être ces données et que, souvent, chaque variable a un poids identique dans l'analyse. L'ajout d'une variable moins significative mais permettant d'expliquer une petite variation, comme c'est possible en économétrie, risquerait ici de fausser le résultat, car elle prendrait une importance comparable aux variables principales de l'analyse. On va donc chercher des données permettant de décrire chacune à leur façon, mais toujours pertinemment, une action. Si on suit M.Aglietta[AGL 98], il y a différents types d'acteur sur les marchés, qui correspondent à des comportements basés sur des analyses différentes : l'analyse technique, l'analyse fondamentale, le suivi de tendance. Chacune de ces analyses correspond à un comportement d'acteur qui agit sur l'évolution du titre et est donc pertinente pour décrire ce titre.

- Analyse technique : Pratiquée par les spéculateurs professionnels avec une vision de court terme ou de très court terme. Comme ils recherchent et obtiennent des informations fines sur les mouvements des cours, ils tiennent compte des impulsions données au marché par les autres groupes d'opérateurs, dans la mesure où ces impulsions affectent le prix. Les données pertinentes à rechercher sont donc les évolutions des cours, de la façon la plus fine possible, l'idéal étant de disposer des courbes intraday du cours et du volume des échanges, voir du carnet d'ordre pour des moments précis (clôture ou ouverture).

- Analyse fondamentale : Pratiquée par les investisseurs institutionnels gérant de très grosses masses d'épargne contractuelle sur le moyen ou long terme. Elle est basée sur l'idée que les marchés sont au moins faiblement efficients à moyen terme. Les prix des actions sont une représentation de l'espérance de gain futur et du risque lié à la volatilité de l'action. Les données pertinentes sont donc liées ici aux perspectives de croissance de l'entreprise et à la volatilité de l'action. La capitalisation de la société, facteur de stabilité, peut être prise en compte. De même que le secteur d'activité, qui peut avoir plus ou moins de perspectives de développement. Le mieux est encore d'utiliser des données issues d'une analyse financière des comptes sociaux sur plusieurs années. Ces comptes sociaux sont publics, et souvent disponibles sur les sites des entreprises. Enfin, la volatilité du capital de l'entreprise de même que l'écart type des performances de l'action représentent de bons estimateurs du risque lié au titre.

- Suivi de marché : Une bonne part des investisseurs institutionnels ou des particuliers qui ont des fonds à placer n'est pas informée. Ils tentent donc de suivre la tendance et contribuent ainsi à l'amplifier. La prise en compte de ce comportement peut se faire en intégrant les performances passées du titre dans l'analyse.

### 3.2.2. Recherche des données et intégration à Excel

La plupart des données citées précédemment sont disponibles ou au moins calculables à partir de données disponibles gratuitement sur Internet. Le problème consiste ensuite à les intégrer sous une forme utilisable par SODAS, c'est à dire concrètement les obtenir sous la forme d'une base de données.

Les données techniques (cours et volume) sont les plus faciles à obtenir car disponibles directement sous forme de table à télécharger. Pour un simple problème de taille, l'analyse est limitée ici à 280 valeurs sur 5 mois, soit 104 jours boursiers. Les cours (ouverture, min, max, clôture, volume) ont été téléchargés sous format texte mois par mois sur le site des Echos. Pour être ensuite utilisables, ces données sont ensuite rassemblées dans une même table à l'aide d'une macro.

Les données fondamentales concernant les actions, telles que le nombre total de titres, le secteur, les poids dans les indices, sont plus difficiles à obtenir sous forme de base de données. Elles sont très souvent affichées, mais à l'intérieur de pages décrivant le titre, et non dans un listing avec tous les autres titres. Le seul moyen pour les récupérer est donc de copier chacune de ces pages récapitulant les informations sur un titre dans une table Excel pour ensuite pouvoir les récupérer automatiquement grâce à une macro, en profitant du fait que toutes ces pages ont le même format. Chaque site boursier propose une page récapitulant les informations sur chaque titre. Le plus complet (notamment parce qu'il offre les trois niveaux de secteur) et le plus facile à traiter sous Excel est celui de la bourse de Paris. Les 280 pages du site décrivant les actions ont donc été collées sur une feuille Excel et trois macros sont ensuite allées chercher sur les 280 feuilles les informations concernant le secteur, le marché et les données décrivant le titre.

### 3.2.3. Exportation des données SODAS à travers DB2SO et Access

Une fois les données disponibles sous forme de tables Excel, il faut les lier sous Access et les structurer à l'aide de requêtes pour pouvoir ensuite les importer sous DB2SO qui les convertira au format SODAS. DB2SO nécessite des requêtes différentes selon le type de variables à importer : une requête pour les variables intervalles, une pour les variables multivaluées et une pour les variables simples. Pour chaque type d'objet symbolique (action, secteur, …) on construit donc entre une et trois requêtes (selon qu'on a ou non tous les types de variables), et on l'importe dans DB2SO. On construit également une requête pour la taxonomie des secteurs en trois niveaux, qu'on peut appliquer pour les objets symboliques de type action. Le fichier DB2SO peut alors être exporté vers SODAS.

### 3.3. *Analyse des données à l'aide de SODAS*

Une fois le fichier SODAS souhaité disponible, on applique les méthodes permettant de répondre au type de description qu'on souhaite.

Par exemple, si on veut obtenir une description générale du secteur informatique depuis quatre mois, on commence par créer le fichier SODAS correspondant (on reprend les requêtes Access en spécifiant « INFORMATIQUE » comme paramètre de secteur, on importe les requêtes dans DB2SO et on exporte le fichier SODAS), puis on peut par exemple appliquer la méthode DIV en ne sélectionnant que les paramètres fondamentaux et de long terme (les paramètres de court terme troubleraient l'analyse), qui permettra de diviser le secteur en groupes d'actions ayant des caractéristiques similaires, permettant de mieux comprendre la structure globale du secteur.

Les différentes méthodes de SODAS répondent en fait à des questions différentes ; par exemple, appliquées à l'analyse d'un secteur en prenant les actions comme objet symbolique et en retenant les variables fondamentales et de long terme:

- DIV, qui renvoie une classification hiérarchique, permet de répondre à la demande « donne-moi les groupes d'actions principaux caractérisant ce secteur ainsi que les critères de sélection »
- PCM, qui réalise une analyse en composante principale, répond à « affiche-moi l'ensemble des actions du secteur de façon à faire ressortir les différences »
- FDA, qui réalise une analyse discriminante, répond à « affiche-moi l'ensemble des actions de façon faire ressortir les groupes qui s'oppose selon tel critère, par exemple la performance sur un mois »
- TREE, qui réalise un arbre de décision, répond à « donne-moi les critères qui permettent de classer au mieux les actions en fonction de telle variable, par exemple la performance sur un mois »
- PYR, qui réalise une pyramide permet de répondre à « affiche-moi toutes les actions du secteur et trace la pyramide de façon à faire apparaître les niveaux de dissimilarité entre chaque couple d'action »

### 4. Résultats

Du fait du très grand nombre de combinaisons niveau d'analyse - objet symbolique - variables – méthode qu'il est possible de choisir, SODAS permet de répondre à un nombre quasi-illimité de question, ce qui renforce encore l'importance du choix et de la pertinence de ce choix pour l'utilisateur. Le grand nombre de cas possible empêche également une présentation exhaustive des possibilités et on va se contenter ici de présenter deux exemples significatifs d'application qui ont étés réalisés sur 250 valeurs entre le 1$^{er}$ novembre 1999 et le 27 mars 2000.

### 4.1. *Analyse globale*

La première analyse est la plus générale et la plus simple : « Décris-moi globalement le marché depuis quatre mois et plus particulièrement les événements récents». Concrètement, on se place au niveau du marché global, on considère des objets symboliques décrivant des actions, on utilise les variables fondamentales et de moyen-long terme et on applique la méthode DIV pour obtenir, par exemple, une classification en 8 classes (on utilise une normalisation en fonction de l'inverse de l'écart type). Le résultat obtenu est :

```
PARTITION IN 8 CLUSTERS :
------------------------:

Explicated inertia : 69.109044

                         +----  Classe 1 (Ng=2)
                         !
              !----6- [perfmois <= -51.945099]
              !    !
              !    +----  Classe 7 (Nd=140)
              !
         !----4- [perf2sem <= 5.065145]
         !    !
         !    +----  Classe 5 (Nd=53)
         !
    !----2- [volat20 <= 4.025935]
    !    !
    !    !              +----  Classe 3 (Ng=21)
    !    !              !
    !    !    !---7- [perf2sem <= -4.492455]
    !    !    !         !
    !    !    !         +----  Classe 8 (Nd=18)
    !    !    !
    !    !----3- [volat10 <= 25.455649]
    !         !
    !         +----  Classe 4 (Nd=1)
    !
!----1- [capim10 <= 90903948.000000]
    !
    !    +----  Classe 2 (Ng=14)
    !    !
    !----5- [capitmds <= 150.324997]
         !
         +----  Classe 6 (Nd=1)
```

**Figure 1.** *Résultat de l'application de la méthode DIV sur 250 actions*

La première réaction en voyant le résultat obtenu lors de la classification est que le programme a trouvé des classes évidentes, et c'est précisément ce qui est remarquable. En effet, chacune des divisions correspond à une caractéristique bien connue du marché, soit constante, soit lors des derniers mois. C'est la preuve que ces caractéristiques constituent des *informations de synthèse pertinentes* qui n'avaient pas été fournies au programme et qu'il a donc retrouvé tout seul. Il a de plus fourni les plus évidentes, ce qui signifie aussi celles qui caractérisent le mieux le marché :

- Capitaux moyens échangés par jour sur deux semaines < 90 M€ (1). Cette classification permet de séparer les 15 plus grosses valeurs de la cote qui ont un comportement nettement différent des autres notamment du fait de leur forte liquidité et de leur poids dans l'indice. Cette division, qui explique à elle seule 26% de l'inertie totale correspond effectivement à la partition fondamentale qu'aurait donnée n'importe quel spécialiste des marchés pour décrire le comportement des valeurs du SBF250. On peut d'ailleurs remarquer que le critère choisi (les capitaux moyens échangés par jour) est sans doute plus significatif que la capitalisation totale qui vient plus naturellement à l'esprit.

- Capitalisation < 150 milliards d'euros (5). Parmi les 15 plus grosses valeurs de la cote, le programme sépare seule France Télécom. Cette action permet d'expliquer à elle seule 6% de l'inertie totale. On retrouve le phénomène de la taille. Le poids de France Télécom, aussi bien en terme de capitalisation totale, que de capitaux échangés ou de poids dans l'indice, conduit l'action à se comporter de façon très particulière. Ce phénomène a été encore renforcé au cours du mois de mars par l'annonce de l'introduction en bourse de sa branche Internet Wanadoo, ce qui a conduit l'action à augmenter de plus de 15% dans une journée, phénomène très inhabituel pour une société de cette taille.

- Volatilité moyenne du capital sur 1 mois < 4‰ (2). Cette division, qui explique 14% de l'inertie, sépare les quarante entreprises dont le capital est le plus volatil des 195 restantes. La forte volatilité du capital conduit à un comportement différent : plus fortes variations du titre, mouvements spéculatifs, possibles changement de propriétaires, …

- Volatilité moyenne du capital sur 2 semaines < 25‰ (3). SODAS choisit ici d'isoler Moulinex dans une classe à part. Ce titre explique à lui seul plus de 8% de l'inertie, plus que France Télécom (qui sera détaché à l'étape 5). Il s'agit là encore d'une information tout à fait pertinente et caractérisant bien le marché au cours du dernier mois. Moulinex a en effet rencontré des problèmes sociaux qui ont fait fortement chuter le cours de l'action. Puis les dirigeants ont annoncé une augmentation de capital pour fin mars. Parallèlement, un groupe italien a profité de la forte volatilité du capital pour accroître sa part dans le groupe jusqu'à atteindre 25% et souhaitait fin mars fusionner Moulinex avec sa filiale Brandt. La remarque faite par SODAS sur la forte volatilité de Moulinex est donc tout a fait justifiée et constitue un des faits majeur du marché au cours des mois étudiés.

- Performance moyenne sur 2 semaines < -5%. Ce critère permet au programme d'isoler, dans les 38 valeurs à forte volatilité, les 21 valeurs ayant subi la chute des valeurs technologique au cours des deux dernières semaines. On retrouve ainsi toutes les valeurs informatiques qui, après avoir fortement augmenté de concert, ont subi la correction de mars. Au contraire, les 18 restantes, tout en gardant une volatilité élevée, n'ont pas subi la chute des marchés technologiques.

- Performance moyenne sur 2 semaines < 5%. Les 195 actions restantes sont divisées en deux groupes de 142 et 53 actions, trop importants pour pouvoir en tirer des conclusions.

- Performance moyenne sur 1 mois < -52%. Beneteau et Europstat sont isolés du reste du groupe. En réalité, cela provient d'une erreur dans les données. Les deux titres ont fait l'objet d'une division du nominal, ce qui a été interprété par le programme comme une forte chute des cours.

Ces caractéristiques, souvent évidentes, sont en fait une confirmation du bon fonctionnement de l'analyse : l'analyse donnée par le programme correspond à celle qu'aurait pu donner une personne connaissant correctement le fonctionnement de la bourse et ayant suivi l'actualité du dernier mois. Cela signifie également que sur des secteurs où on ne s'attend pas à ce qu'on va trouver, il n'y a pas de raison que la pertinence des informations obtenues soit moindre, et on va donc obtenir une information nouvelle et de qualité qui nous était, cette fois-ci, inconnue.

**4.2.** *Analyse intersectorielle*

Après avoir décrit le marché en terme d'actions, on peut chercher à le décrire par rapport aux secteurs. Le niveau d'analyse reste le marché global. Les objets symboliques décrivent alors des secteurs (il y a 22 objets symboliques, ou 12 si on se place au niveau 2). Les variables considérées sont les même que précédemment mais agrégées sous forme de variables intervalles pour chaque secteur. Pour obtenir une représentation graphique des secteurs sur le marché, on peut justement utiliser l'outil d'analyse en composante principale PCM, en lui fournissant les données intervalle concernant les données fondamentales et de moyen terme. On obtient le graphique suivant avec 22 secteurs et en considérant les 2 premiers axes (qui expliquent respectivement 20 et 19% de l'inertie)

**Figure 2.** *Résultat de l'analyse factorielle sur 22 secteurs, axes 1 et 2*

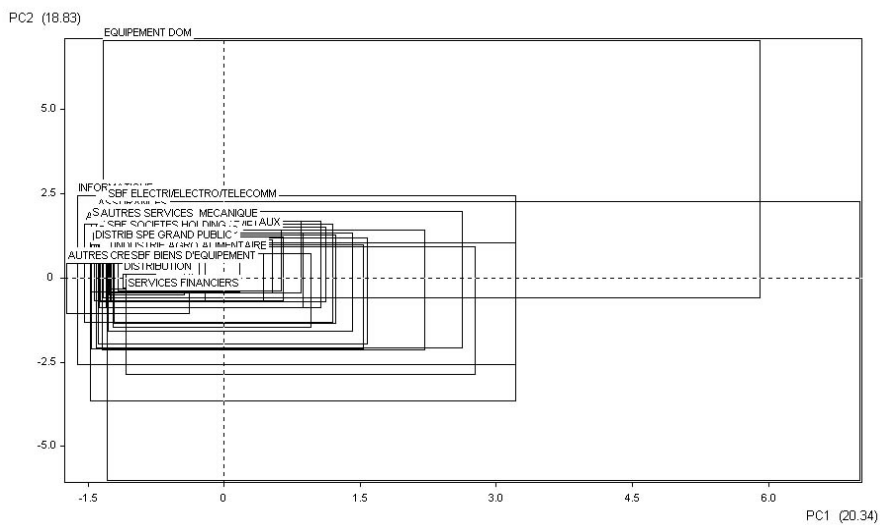

On s'aperçoit ici que les secteurs sont tous, au moins en partie, regroupés autour de l'origine. Aucun secteur ne se détache ou s'oppose à un autre. Les comportements des secteurs ont tendance à se rapprocher d'un « comportement standard » qu'on peut rapprocher du « portefeuille de marché » de la théorie financière. Le comportement d'un secteur entier ne peut s'éloigner de celui du marché, seules certaines valeurs s'en éloignent temporairement et étirent la zone couverte par leur secteur, comme France Télécom avec les télécommunications ou Moulinex avec les équipements domestiques. Cet étirement est accentué par le fait qu'on considère des valeurs intervalle, et non des lois qui auraient fait ressortir la particularité de ces cas. En fait, ce n'est pas le secteur des télécommunications qui est différent de celui des équipements domestiques, c'est France Télécom qui s'est comporté très différemment de Moulinex. L'intérêt de ce graphique n'est donc pas dans les oppositions mais dans leur absence. Certains secteurs apparaissent en effet très homogènes alors que d'autres apparaissent plus dispersés. On peut l'apercevoir de façon plus nette si on se limite à 12 secteurs. De plus, on peut annuler les effets France Télécom et Moulinex en se basant sur les axes 3 et 4 :

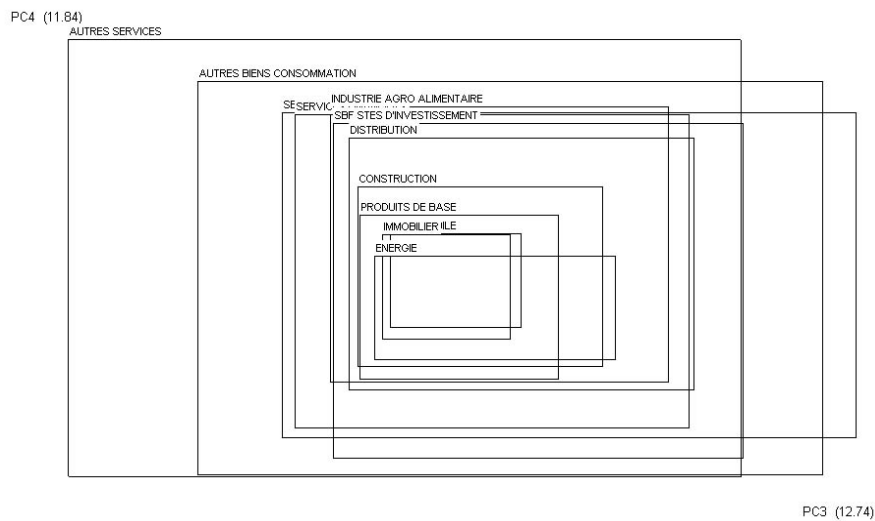

**Figure 3.** *Résultat de l'analyse factorielle sur 12 secteurs, axes 3 et 4*

Les autres services, comprenant notamment l'informatique, apparaissent comme le secteur le plus complexe et comportant tous les comportements possibles. Au contraire, des secteurs comme l'énergie, l'automobile, l'immobilier et les produits de base sont non seulement très homogènes en interne mais également entre eux, ils se comportent tous de façon similaire. On peut avoir confirmation de cette analyse en traçant la pyramide correspondante :

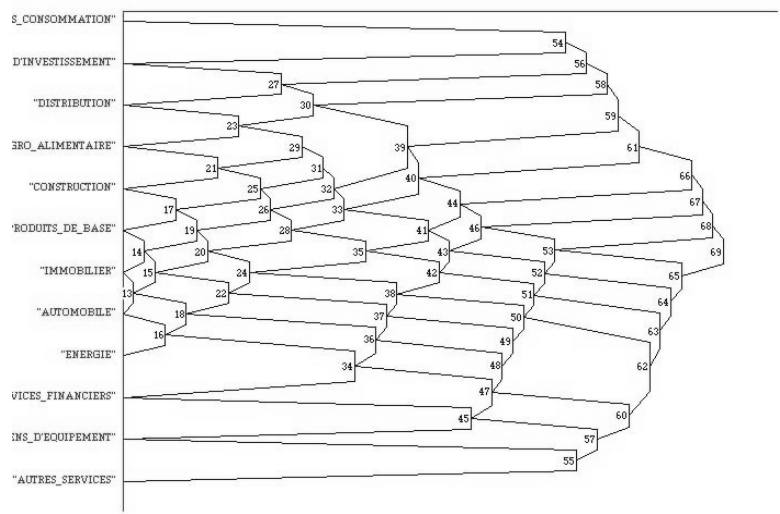

**Figure 4.** Résultat de l'application de la méthode PYR sur 12 secteurs

On retrouve les secteurs de l'immobilier, de l'automobile, des produits de base et de l'énergie comme très proche. Au contraire, les biens de consommation (Moulinex), les autres services (informatique) et les biens d'équipement (France Telecom) sont non seulement très différents des secteurs « du centre » mais également très différents entre eux : le palier 45 qui est le premier palier liant un de ces trois secteurs à un autre, est juste en dessous du palier décrivant l'écart entre les deux secteurs les plus éloignés parmi les 9 restant (le 46).

## 5. Perspectives

Les tests qui ont été fait ici montrent que les données boursières se prêtent bien à l'analyse de données symboliques. La principale contrainte devient, pour chaque nouvelle question, de devoir réaliser les requêtes correspondantes dans Access (changer le paramètre), importer les résultats dans DB2SO et exécuter SODAS. Il faut de plus remettre à jour régulièrement les données des tables Excel. Une évolution logique de cette application serait donc d'automatiser tous ces processus pour obtenir un programme ou une page Web simple à utiliser.

Une autre direction pour l'amélioration de cette analyse serait l'intégration de plus de données financières. On n'a utilisé ici que des données boursières, alors que des données financières telles que les taux de rendement, d'endettement, de rentabilité,… sont autant d'informations importantes dans l'évaluation des sociétés et donc dans l'explication et l'évolution des cours.

Enfin, on peut imaginer les perspectives intéressantes qu'offriraient certaines extensions de l'analyse de données symboliques qui ne sont pas encore dans SODAS. On peut par exemple penser à la régression symbolique, qui permettrait non plus de décrire, mais de prévoir l'évolution des cours grâce aux données symboliques…

## 6. Bibliographie